\documentclass [12pt] {article}
\usepackage{amsfonts}
\usepackage{amssymb}

\newcommand{\qed} {\hspace {0.1in} \rule {1.5mm} {3.5mm}}

\newtheorem{theorem}{Theorem}

\newtheorem{proposition}{Proposition}[section]

\def\lnt{\lim_{n\to\infty}}

\def\limi{\lim_{i\to\infty}}
\def\limn{\lim_{n\to\infty}}
\def\gseq{\{G_n\}_{n=1}^\infty}

\def\ess{\mbox{Ess}\,\mbox{sup}\,}

\def\deg{\mbox{deg}\,}

\def\tr{\mbox{Tr}\,}
\def\mat{\mbox{Mat}\,}
\def\trr{\tr_{\cG}}

\def\<{\langle}
\def\>{\rangle}

\def\proof{\smallskip\noindent{\it Proof.} }

\def\bN{{\mathbb N}}
\def\bC{{\mathbb C}}
\def\deg{\mbox{deg}\,}

\def\cA{\mbox{$\cal A$}}
\def\cB{\mbox{$\cal B$}}
\def\cC{\mbox{$\cal C$}}
\def\cH{\mbox{$\cal H$}}
\def\cG{\mbox{$\cal G$}}

\def\crd{\cC^{r,d}}

\def\cR{\mbox{$\cal R$}}

\def\cN{\mbox{$\cal N$}}

\def\to{\rightarrow}

\begin{document}
\title{Weak convergence of finite graphs, integrated density of states
 and a Cheeger type inequality} 
\author{\sc G\'abor Elek
\footnote {The Alfred Renyi Mathematical Institute of
the Hungarian Academy of Sciences, P.O. Box 127, H-1364 Budapest, Hungary.
email:elek@renyi.hu, Supported by OTKA Grants T 049841 and T 037 846}}
\date{}
\maketitle
\vskip 0.2in
\noindent{\bf Abstract.}  In \cite{Elek} we proved that the limit
of a weakly convergent sequence of finite graphs can be viewed
as a graphing or a continuous field of infinite graphs. Thus one can
associate a type $II_1$-von Neumann algebra to such graph sequences.
We show that in this case the integrated density of states exists that
is the weak limit of the spectra of the graph Laplacians of the finite
graphs is the KNS-spectral measure of the graph Laplacian of
the limit graphing. Using this limit technique we prove a Cheeger type
inequality for finite graphs.

\vskip 0.2in
\noindent{\bf AMS Subject Classifications:} 05C80, 46L10
\vskip 0.2in
\noindent{\bf Keywords:}  weak convergence of graphs,
von Neumann algebras, isoperimetric
inequalities, integrated density of states
\vskip 0.3in
\newpage
\section{Introduction}
\subsection{Weak convergence and limits of colored graph sequences}
First let us recall some of the definitions and the main result from
$\cite{Elek}$. A {\it rooted colored $d$-graph} is a finite simple graph $G$
with the following properties :
\begin{itemize}
\item
$G$ has a distinguished vertex (the root).
\item
For any $p\in V(G)$, $deg(p)\leq d$.
\item
The edges of $G$ are properly colored by the colors $c_1, c_2,\dots, c_{d+1}$.
That is for each vertex $q\in V(G)$ the outgoing edges from $q$
are colored differently. Note that by Vising's theorem each graph with
vertex degree bound $d$ has such an edge-coloring.
\end{itemize}
A {\it rooted $(r,d)$-ball} is a rooted colored $d$-graph
such that $\sup_{y\in V(G)} d_G(x,y)\leq r$, where $x$ is the root of
$G$ and $d_G$ denotes the shortest path distance. Two rooted colored
$d$-graphs
$G$ and $H$ are rooted isomorphic if there exists a graph isomorphism
between them preserving the colors and mapping one root to the other one.
We denote by $\crd$ the finite set of rooted isomorphism classes of
rooted $(r,d)$-balls. Now if $G$ is a rooted colored $d$-graph and $\cA\in
\crd$, then $T(G,\cA)$ denotes the set of vertices $v$ in $V(G)$ such
that $\cA$ represents the rooted isomorphism class of the
$r$-neighborhood of $v$, $B_r(v)$.  Set $p_G(\cA):=\frac{T(G,\cA)}{|V(G)|}\,.$
That is $G$ determines a probability distribution on $\crd$, for any $r\geq
1$.
Now let $\{G_n\}_{n=1}^\infty$ be a sequence of finite simple connected
graphs with vertex degrees bounded by $d$, edge-colored by the colors
$c_1, c_2,\dots, c_{d+1}$. Suppose that $V(G_n)\to\infty$. Following
Benjamini and Schramm \cite{BS} we say that $\{G_n\}_{n=1}^\infty$
is {\it weakly convergent}
if for any $r\geq 1$ and $\cA\in\crd$, $\lnt p_{G_n}(\cA)$ exists.

\noindent
Now let $X$ be a compact metric space with a probability measure $\mu$.
Suppose that $T_1,T_2,\dots, T_{d+1}$ are measure preserving continuous
involutions on $X$ satisfying the following conditions :
\begin{itemize}
\item For any $p\in X$, $T_i(p)=p$, for at least one $i$.
\item If $T_i(p)=T_j(p)$, $i\neq j$, then $T_i(p)=T_j(p)=p$.

\end{itemize}
We call such a system $\cG=\{X, T_1,T_2,\dots T_{d+1},\mu \}$
a {\it $d$-graphing}. A  $d$-graphing determines an equivalence relation on
the points of $X$. Simply, $x\sim_{\cG} y$ if there exists a sequence of
points
$\{x_1, x_2,\dots, x_m\}\subset X$ such that
\begin{itemize}
\item $x_1=x, x_m=y$
\item $x_{i+1}=T_j(x_i)$ for some $1\leq j \leq d+1$.
\end{itemize}
Thus there exist natural simple (generally infinite)
 graph structures on the equivalence classes,
the {\it leafgraphs}. Here $x$ is adjacent to $y$, if $x\neq y$ and
$T_j(x)=y$.
By our conditions, all the leafgraphs have vertex degrees bounded by $d$
and are naturally edge-colored (the $(x,y)$ edge is colored by
$c_j$ if $T_j(x)=y$).
Now if $\cA\in\crd$, we denote by $p_{\cG}(\cA)$
the $\mu$-measure of the points $x$ in $X$ such that the rooted
$r$-neighborhood of $x$ in $\cA$ in its leafgraph is isomorphic to $\cA$
as rooted $(r,d)$-ball.
We say that $\cG$ is the limit graphing of the weakly convergent sequence
$\gseq$, if for any $r\geq 1$ and $\cA\in\crd$
$$ \limn p_{G_n}(\cA)= p_{\cG} (\cA)\,.$$
In \cite{Elek} we proved the following result.
\begin{theorem}
\label{t1}
If $\gseq$ is a weakly convergent system of finite connected $d$-graphs
as above, then
there exists a $d$-graphing $\cG$ such that for any $r\geq 1$ and
$\cA\in\crd$, $\limn p_{G_n}(\cA)= p_{\cG} (\cA)\,.$
\end{theorem}

\subsection{Cheeger type isoperimetric inequalities}
Now let us recall two basic Cheeger-type isoperimetric inequalities
for graphs. Let $G(V,E)$ be a finite connected graph. For
a finite connected spanned subgraph $A\subseteq G$
we denote by $\partial A$ the set of vertices $x$ in $V(A)$ such that
$x\in V(A)$, but there exists $ y\notin V(A)$ that $x$ and $y$ are adjacent vertices.
The Cheeger constant $h(G)$ is defined as
$$h(G):= \{\inf \frac {|\partial A|} {|A|}\,\mid \,A\subset G,\,|V(A)|\leq
\frac{|V(G)|}{2}\}\,.$$
The classical isoperimetric inequality can be formulated the following
way: For any $\epsilon>0$ and $d\in \bN$, there exists a real constant
$C(\epsilon,d)>0$ such that if $G$ is a finite connected graph with vertex
degree bound $d$ and $\lambda_1(G)\leq C(\epsilon,d)$ then $h(G)\leq
\epsilon$. Here $\lambda_1(G)$ is the first non-zero Laplacian eigenvalue
of $G$.

\noindent
Now let $G(V,E)$ be an infinite connected graph with bounded vertex degrees.
Then the isoperimetric constant of $G$ is defined as follows.
$$i(G):=\{\inf \frac {|\partial A|} {|A|}\,\mid \,A\subset G, \mbox
{$A$ is a finite connected spanned subgraph}\}\,.$$
Then again we have a Cheeger-type inequality (see e.g. \cite{Soa}, Theorem
4.27): For any $\epsilon>0$ and $d\in \bN$, there exists a real constant
$D(\epsilon,d)>0$ such that if $G$ is an infinite connected graph with
vertex degree bound $d$ and 
$$[0,D(\epsilon,d)]\cap \mbox{Spec}\,\Delta_G \neq 0\,,$$
then $i(G)\leq\epsilon$. That is small values in the spectrum indicates
the existence of some subsets with small boundaries. The main result
of our paper can be stated informally the following way:
The abundance of small eigenvalues in a finite graph indicates the
existence of a large system of disjoint subsets with small boundaries.
To formulate our result precisely we need some definitions.
For $\delta>0$ and a finite connected graph $G$
let $S(G,\delta)$ denote the set
of Laplacian eigenvalues (with multipicities) not greater than $\delta$.
Then
$$s(G,\delta):=\frac{|S(G,\delta)|}{|V(G)|}\,.$$
For a finite, connected graph $G$, $\epsilon>0$ and $k\in \bN$ we denote by 
$H(G,\epsilon, k)$ the set of vertices in $G$ which can be covered
by a  connected spanned subgraph $A$ such that $|A|\leq k$ 
and $\frac {|\partial A|} {|A|}\leq
\epsilon$. Then $h(G,\epsilon,k):=\frac {|H(G,\epsilon, k)|}{|V(G)|}.$
Also, we denote by $m(G,\epsilon,k)$ the cardinality of the maximal
disjoint system of  connected spanned subgraphs $A$ in $V(G)$ such that 
$|A|\leq k$ and $\frac {|\partial A|} {|A|}\leq
\epsilon$.
Clearly, there exists a constant $c(\epsilon, d, k)>0$ such that for all
graphs $G$ with vertex degree bound $d$,
$$m(G,\epsilon,k)\geq c(\epsilon, d, k) h(G,\epsilon,k)\,.$$
Using the weak limit technique we prove the following result :
\begin{theorem}
\label{t2}
For any $\epsilon>0$ and $d\in \bN$, there exists a constant $E(\epsilon,d)>0$,
such that the following holds:

\noindent
For any $t>0$, there exist $k(t),n(t)\in \bN$ and $s(t)>0$ such that if
$|V(G)|\geq n(t)$ holds for a finite connected graph $G$ with vertex
degree bound $d$, and $s(G, E(\epsilon,d))\geq t$, then
$m(G,\epsilon, k(t))\geq s(t)\,.$
\end{theorem}
The main idea is to associate a type $II_1$-von Neumann
algebra for a weakly convergent graph sequence and prove that
the integrated density of states exists that is
the spectra of the graph Laplacians of the finite graphs converge
weakly to the so-called KNS-spectrum of the graph Laplacian on
the limit graphing.

\noindent
Note that
in \cite{LS1} Lovasz and Szegedy introduced the notion of
weak convergence and the limit object for
dense graph sequences. They used the limit technique to give a new proof
of a theorem
of Alon and Shapira in \cite{LS2}. An another application of the limit
technique can be found in the paper of Aldous and Steele \cite{Ald}.

\section{The KNS-measure}

First let us recall the classical definitions of measureable fields
of Hilbert-spaces and operators \cite{Dix}.
Let $Z$ be the standard Borel space with a probability
 measure $\mu$. Let $\cH=l^2(\bN)$
be the complex separable infinite dimensional Hilbert-space with orthonormal
basis $\{e_i\}^\infty_{i=1}$. A measurable field of Hilbert-spaces is
given by a sequence $\{f_n\}^\infty_{n=1}\subset L_2(Z,\mu)$ such
that $\sum^\infty_{n=1} \int_Z |f_n(x)|^2\,d\mu(x)<\infty\,.$
By the Beppo-Levi Theorem for almost every $x\in Z$:
$f_x=\sum^\infty_{n=1}f_n(x)e_n\in\cH$. We denote the space of the
measurable fields of Hilbert spaces by
$\int_Z\cH_x\,d\mu(x)\,$, which is
 a Hilbert-space with respect to the pointwise
inner product.

\noindent
A field of bounded linear operators is a $\cB(\cH)$-valued
function on $Z$ such that \\ $x\to <A_x(e_i),e_j>$ is a measurable
function for all $i,j$ and $\ess_{x\in Z} \|A_x\|<\infty.$
Such fields of operators form the von Neumann algebra $\int_Z
\cB(\cH)\,d\mu(x)$, with pointwise addition, multiplication and $*$-operation.
Note that $\int_Z A_x\,d\mu(x)\in \int_Z
\cB(\cH)\,d\mu(x)$ is invertible if and only if $A_x^{-1}$ exists
for almost every $x\in Z$ and 
$\ess_{x\in Z} \|A_x^{-1}\|<\infty.$

\noindent
Now we introduce the notion of continuous field of graphs associated
to a $d$-graphing. Let $\cG=(X,T_1,T_2,\dots,T_{d+1},\mu)$ be a $d$-
graphing of our Theorem \ref{t1}.
For each $r\geq 1$, and $\cA_r\in\crd$ we label the vertices of
$\cA_r$ by natural numbers inductively, satisfying the following
conditions:
\begin{itemize}
\item The root is labeled by $1$.
\item Any vertex has different labeling.
\item The labeling of $\cA_r$ is compatible with the labeling of $\cA_{r-1}$.
\item If $d_G(x,y)<d_G(x,z)$, then the label of $y$ is smaller than
the label of $z$.
\item If $|\cA_r|=k$, then the labels of the vertices is exactly
$\{1,2,\dots,k\}$.
\end{itemize}
Therefore for each $x\in X$
we associate an infinite graph with edge-coloring by
$c_1, c_2, \dots, c_{d+1}$ and an extra vertex labeling
by the natural numbers. Of course this graph is isomorphic to the leafgraph
of $x$ and each vertex can also be viewed as an element of $X$.

\noindent
This is the {\it continuous field of infinite graphs}
associated to the graphing.
The total vertex set of this field is $\cR\subset X\times X$, $(x,y)\in\cR$
if $x\sim y$, the equivalence relation given by the graphing. 
 The space
$\cR$ is equipped with the counting measure $\nu$,\cite{FM} hence using the
vertex labels one can view $L_2(\cR,\nu)$ as $\int_X l^2(V(\cG_x)) d\mu(x)$.

\noindent
The leafwise Laplacian operator $\int_X \Delta_{\cG_x}\,d\mu(x)$
is a measurable field of bounded operators on $\int_X l^2(V(G_x)) d\mu(x)$.
Feldman and Moore \cite{FM} introduced an important subalgebra of the
algebra of the fields of bounded operators in the case of graphings.
We briefly review their construction.
Consider the space of bounded
 measurable functions $\kappa:\cR\to \bC$ with 
{\it finite bandwidth}, that is 
 for some constant $b_{\kappa}>0$
depending only on $\kappa$:  $\kappa(x,y)=0$ if
$d_{\cG}(x,y)>b_{\kappa}$. Here $d_{\cG}(x,y)$ denotes the shortest
path-distance on the leaf-graph of $x$. One can associate a measurable
field of bounded operators to $\kappa$ the following way.
 If $f\in L^2(\cR,\nu)$, then:
$$T_{\kappa}(f)(x,y)=\sum_{z\sim x} f(x,z) \kappa (z,y)\,.$$
These operators all called {\it random operators}.
The weak closure of such operators
 $T_{\kappa}$ in $\int_X \cB(l^2(V(\cG_x))d\mu(x)\,$ is the Feldman-Moore
algebra $\cN_{\cR}$. The point is that $\cN_{\cR}$ possesses a trace. It is
a von Neumann algebra of type $II_1$.

The trace of $T_{\kappa}$ is given by
$$\trr(T_{\kappa})=\int_X \kappa(x,x) d\mu(x)\,.$$
The leafwise Laplacian operator on $\cG$, $\Delta_{\cG}$ is clearly
an element of $\cN_{\cR}$, given by a bounded measurable function
of finite bandwidth, where
$$\Delta_{\cG}(f)(x,y)=\deg(y) f(x,y)-\sum_{\mbox{$z$, $z$ is adjacent to $y$}} f(x,z)\,.$$
We shall denote the Laplacians on the finite graphs $G_i$ by $\Delta_{G_i}$.
(Of course the Laplacian does not depend on the edge coloring.)
\begin{proposition}
\label{pu14}
 For any $n\geq 1$:
$$\limi \frac{\tr_i (\Delta^n_{G_i})}{|V(G_i)|}=\trr (\Delta^n_{\cG})\,,$$
where $\tr_i:\mat_{V(G_i)\times V(G_i)} (\bC)\to \bC$ is the usual trace.
\end{proposition}
\proof
Obviously,
$$\tr_i (\Delta_{G_i})=\sum_{x\in V(G_i)} \deg(x)\,$$
and for the powers of the Laplacians:
$$\tr_i (\Delta_{G_i}^n)=\sum_{x\in V(G_i)} s_n(x)\,,$$
where $s_n(x)$ depends only on the isomorphism class of the $n$ ball
around $x$ in $G_i$.
Hence
$$\frac{\tr_i (\Delta_{G_i}^n)}{|V(G_i)|}=\sum_{A\in \cC^{n,d}}
p_{G_i}(\cA)s_n(\cA)\,,$$
where $s_n(\cA)$ is the value of $s_n$ at the root.
On the other hand,
$$ \trr (\Delta^n_{\cG})= \sum_{\cA\in \cC^{n,d}}
p_{\cG}(\cA)s_n(\cA)\,.$$
Hence our proposition follows. \qed

\noindent
By our vertex bound condition $Spec(\Delta_{\cR})$ and $Spec(\Delta_{G_i})$
for all $i$ are contained in some interval $[0,l]$.
Recall that the spectral measure $\lambda_{T}$ of a positive self-adjoint
operator $A$ on the finite dimensional Euclidean space $\bC^n$ is a point
 measure on
 $[0,\infty)$
defined as follows:
$$\lambda(A)=\frac{\#(\mbox{the eigenvalues of $T$ with multiplicities in
$A$})} {n}\,.$$
By the classical (finite dimensional) spectral theorem:
$$\int^l_0 x^n\,d\lambda_{\Delta_{G_i}}(x)=\frac {\tr_i (\Delta^n_{G_i})}
{|V(G_i)|}\,.$$
On the other hand, by the von Neumann's spectral theorem:
$$\int^l_0 x^n\,d\lambda_{\Delta_{\cR}}(x)=\trr(\Delta^n_{\cG})
\,,$$
where $\lambda_{\Delta_{\cG}}[0,x]=\trr (E_{\Delta_{\cG}})([0,x])$
for the projection
valued measure $E_{\Delta_{\cG}}$ associated\\ to the bounded self-adjoint
operator $\Delta_{\cG}$. The measure $\lambda_{\Delta_{\cG}}$ is 
called the KNS-measure
associated to the Laplacian $\Delta_{\cG}$.
 (see \cite{Gri} for a discussion on spectral
measures.)

\noindent
That is for any polynomial $p\in\bC[x]$,
$$\limi \int^l_0 p(x) d\lambda_{\Delta_{G_i}}(x)= \int^l_0 p(x) 
d\lambda_{\Delta_{\cR}}(x)\,.$$
Therefore we proved that the integrated density of states exists. (see
\cite{Lenz} for a discussion on integrated density of states and random
operators)
\begin{theorem}\label{t3}
If the sequence of graphs $\{G_i\}_{i=1}^\infty$ weakly converges
to the graphing $\cG$, then the associated spectral measures 
$\lambda_{\Delta_{G_i}}$ weakly converge to $\lambda_{\Delta_{\cR}}$.
(see also \cite{Ser})
\end{theorem}

\section{The proof of Theorem \ref{t2}}

Let us choose $E(\epsilon,d):=\frac{1}{2} D(\epsilon,d)$, where
$D(\epsilon,d)$ is the constant in the isoperimetric inequality for
infinite graphs.
Suppose that the theorem does not hold. Then there exists $t>0$ with
the following property:
For any $k$, one can choose a sequence $\{G_n^k\}_{n=1}^\infty$ of
finite graphs with vertex degrees bounded by $d$ such that
\begin{itemize}
\item $|V(G_n^k)|\to\infty$\,.
\item $s(G_n^k, \frac{1}{2} D(\epsilon,d))\geq t.$
\item $\limn m(G^k_n,\epsilon,k)=0\,.$
\end{itemize}
Since if $l>k$ then $m(G,\epsilon,k)\leq m(G,\epsilon,l)$, we can choose
a weakly convergent
 sequence of finite connected colored d-graphs $\{H_n\}^\infty_{n\geq 1}$
such that
\begin{itemize}
\item $|V(H_n)|\to\infty$\,.
\item $s(H_n, \frac{1}{2} D(\epsilon,d))\geq t.$
\item $\limn m(H_n,\epsilon,k)=0\,$ for any $k$\,.
\end{itemize}
By the argument in the Introduction one can see that
$\limn h(H_n,\epsilon,k)=0\,$ as well.

\noindent
Therefore if $L$ is a fixed finite connected graph and $\cA_r$ 
($\mbox{diam}\,(L)<r$) contains
a spanned subgraph $L'$ isomorphic to $L$,
 $\frac{|\partial L'|}{|L'|}\leq\epsilon$
such that the root is in $L'$, then $\limn p_{H_n}(\cA_r)=0$.
Consequently, if one consider the limit graphing $\cG$
of Theorem \ref{t1}, then for each leafgraph $\cG_x$,
$i(\cG_x)\geq \epsilon$. That is for each leafgraph
the spectral gap of the Laplacian at the zero is greater than
$\frac{3}{4} D(\epsilon,d)$. Therefore,  the measurable field
of operators $(\int_X \Delta_{\cG_x} \,d\mu(x)-\lambda)$ is invertible,
with a uniformly bounded inverse,
whenever $\lambda\leq \frac{1}{2} D(\epsilon,d)$. In other words
$$\mbox{Spec}(\int_X \Delta_{\cG_x} \,d\mu(x))\cap 
[0,\frac{1}{2} D(\epsilon,d)]=\emptyset\,.$$
Now we use our Theorem \ref{t3}. The spectrum of
$\int_X \Delta_{\cG_x} \,d\mu(x)$ is the same in the von Neumann algebra
$\int_X \cB(l^2(\cG_x))d\mu(x)$ 
as in the Feldman-Moore subalgebra.
However by Theorem \ref{t3}, the KNS-measure of the interval
 $[0,\frac{1}{2} D(\epsilon,d)]$ is at least $t$. This leads to
a contradiction. \qed

\end{document}